\documentclass[a4paper,11pt, twoside]{article}
\usepackage{amsfonts,amsmath,amssymb}
\usepackage{amsthm}
\usepackage{geometry}
\usepackage{graphicx}
\usepackage{hyperref}
\usepackage{breakcites}
\usepackage{booktabs}
\usepackage{url}
\usepackage{tabularx}
\usepackage{color}
\renewcommand{\tabcolsep}{15pt}
\usepackage{vmargin}
\setmarginsrb{2.6cm}{2cm}{2.6cm}{3cm}{0mm}{5mm}{0mm}{10mm}
\usepackage[all]{nowidow}
\newcolumntype{R}{>{\raggedleft\arraybackslash}X}
\usepackage{lastpage}
\usepackage{fancyhdr}
\fancypagestyle{title}{
  \setlength{\headheight}{0pt}%
  \fancyhf{}
  \fancyfoot[L]{\scriptsize \sffamily {\bfseries CONTACT} \; N. Karjanto \; \url{natanael@skku.edu}\\
    \copyright 2017 Informa UK Limited, trading as Taylor \& Francis Group}
  \fancyhead[L]{\scriptsize \sffamily INTERNATIONAL JOURNAL OF MATHEMATICAL EDUCATION IN SCIENCE AND TECHNOLOGY, 2017\\
  VOL. 48, NO. 6, \thepage--\pageref{LastPage}\\ \url{http://dx.doi.org/10.1080/0020739X.2017.1285060}}
}

\begin{document}
\title{\textbf{Attitude toward mathematics among the students at Nazarbayev University Foundation Year Programme}}
\author{\normalsize N. Karjanto\\ 
{\small Department of Mathematics, University College, Natural Science Campus} \\
{\small Sungkyunkwan University, Suwon 16149, Republic of Korea}}
\date{\scriptsize Article History: Received 24 August 2016, Accepted 17 January 2017, Published online 06 March 2017}
\maketitle
\thispagestyle{title}

\begin{abstract}
This article investigates the attitude toward mathematics among the students enrolled in the Foundation Year Programme at Nazarbayev University. 
The study is conducted quantitatively and an inventory developed by Tapia and Marsh II is adopted in this research.
The inventory consists of 40 statements on the five-point Likert scale.
Gender, specialization and final high school score in mathematics are collected. The number of valid returned questionnaires is 108.
There are 55 males, 53 females, 73 Mathematical--Physics (MP), 22 Biology--Chemistry (BC) and 13 International Relations--Economics (IRE) students. 
Generally, they have a positive attitude toward mathematics, with the score mean and standard deviation are 3.999 and 0.531 out of five, respectively.
We confirm a hypothesis on a positive correlation between previous high achievement in mathematics and favorable attitude toward it.
The correlation value is $r = 0.300$, its effect size is medium and it is extremely significant ($p$-value $= 0.0008 < 0.001$).
There is no significant difference between male and female students in terms of their attitude toward mathematics ($t = 0.6804$, $p$-value $= 0.2489 > 0.05$).
There is a very significant difference between students who specialize in IRE and MP in terms of their attitude toward mathematics ($F(2,105) = 5.6848$, $p$-value $= 0.0045 < 0.01$).
\end{abstract}

\section{Introduction}

A success of students in their study often relates to the attitude of the students themselves to the subjects they are studying during their school period,
which spans from elementary school throughout the college and university.
At the tertiary level, the students who intend to major in Natural Science, Computer Science, Engineering and Economics need a solid background in mathematics.
This solid background can be cultivated by internal and external factors, where the former one includes the students' attitude toward mathematics
and a strong motivation to learn the subject further, to dig deeper into the material so to speak.
In this paper, we are interested in investigating the attitude toward mathematics among the students enrolled in the Nazarbayev University Foundation Year Programme
(NUFYP). The study is conducted quantitatively and an inventory developed by~\cite{tapia2004instrument} is adopted and utilized in this research.

Since the Fall semester of 2010, Nazarbayev University (NU) has accepted at least 500 students annually out of over 3,500 applicants into its Foundation program NUFYP. They are the brightest Kazakhstan's secondary school graduates (K-11) and are competitively selected through a three-stage admission process. An intensive academic one-year program at NUFYP prepares them for successful undergraduate study in three NU Schools: Engineering (SEng), Science and Technology (SST) or Humanities and Social Sciences (SHSS). Students' critical thinking and English language skills are developed. They are also encouraged to work independently in order to be well-equipped to enroll in the undergraduate program at NU. All classes are taught in English at both NUFYP and at the undergraduate levels.

The program was formerly known as the Centre for Preparatory Studies (CPS) and had existed for five years from the academic year 2010/2011 until the academic year 2014/2015. The strategic partner for the CPS program was the University College London (UCL) and thus the curriculum adopts the British tertiary education system.
Starting from the academic year of 2015/2016, the curriculum in the Foundation program has been broadened so that all students will experience a liberal arts education with a range of elective subjects in the Humanities and Sciences. Thus, it is renamed as NUFYP. The student body will be expanded to approximately 700 students annually. All students are expected to be proficient in Mathematics and English upon the completion of the program. The new strategic partner is the University of Warwick from the UK. Even though the program still adopts the British education system, the curriculum has been tailored according to the needs of students in Kazakhstan. The University of Warwick is providing an input into the guidance on the learning process as well as the curriculum development.

Since the majority of NU students would enroll in mathematics-related majors (Natural Science, Computer Science, Engineering and Economics), 
it is essential to find out their attitude toward mathematics at the early stage of their study, namely at the Foundation level. 
The exceptions are the students who major in Political Science and International Relation (PSIR) and World Language and Literature (WLL). 
Even for both PSIR and WLL majors, the curricula still require the students enrolled in these majors to fulfill a certain number of credits from mathematics and 
mathematics-related courses. 

There is an abundant literature on students' attitude toward mathematics and the list mentioned in this paper is far from exhaustive.
A classic literature review on attitude toward mathematics has been summarized by~\cite{aiken1970attitudes}.
The author covers several major topics, including methods of measuring attitudes toward arithmetic and mathematics,
the distribution and stability of mathematics attitudes, the effects of attitudes on achievement in mathematics,
the relationships of mathematics attitudes to ability and personality factors, to parental attitudes and expectations, to peer attitudes, and to teacher
characteristics, attitudes and behavior. \cite{aiken1974two} also proposed on two scales of attitude toward mathematics,
where there is a distinction between enjoyment of mathematics and the value of mathematics, where the latter one admits the importance and the relevance
of mathematics to individual and society.

The focus of our study and the corresponding relevant literature are on the Attitudes Toward Mathematics Inventory (ATMI) promoted by~\cite{tapia1996attitudes}.
The inventory is developed to measure students' attitude toward mathematics and to seek the underlying dimensions that comprise it.
The author collected a sample of 544 students taking mathematics at a private bilingual, American curriculum-based,
preparatory school in Mexico City. The instrument initially consists of 49 items and the reliability Cronbach-$\alpha$ coefficient is $\alpha = 0.96$.
After dropping the nine weakest items, the reliability was then increased to $\alpha = 0.97$.
Four factors captured in this instrument are the students' sense of security, the value of mathematics, the motivation to study mathematics and the enjoyment of mathematics. A similar study conducted with 262 middle school students at a private,bilingual college preparatory school in Mexico City has been presented by~\cite{tapia2000attitudes}. The authors concluded a maximum likelihood factor analysis with a varimax rotation providing three factors, namely self-confidence,
the enjoyment of mathematics and the value of mathematics. 

Furthermore, a continuing study has confirmed the results that the four-factor model in the ATMI also holds for college students studying in American universities~\cite{tapia2002confirmatory}. The psychometric properties of the ATMI developed by~\cite{tapia2004instrument} has been examined by~\cite{afari2013examining}. To validate the questionnaire, the author collected data from 269 middle school students in the United Arab Emirates (UAE). The result has provided an additional support to the factorial validity of the ATMI where a high-reliability coefficient for overall inventory was confirmed by factor analysis.

The attitudes toward mathematics of preservice elementary teachers entering an Introductory Mathematics Methods course is examined in a PhD thesis by~\cite{schackow2005examining}. The study sought to determine the extent to which preservice teachers' attitudes toward mathematics changed during the methods course and the correlation between preservice teachers’ initial attitudes toward mathematics and their academic achievement in the methods course. Another PhD dissertation examined the attitudes toward mathematics and their academic outcomes of Beginning Algebra students at a Central Florida Community College~\cite{sisson2011examining}. The author discovered that the students' overall attitudes had a positive change over the semester as indicated by the improvement of their score. A statistically significant association was found in change scores in the ATMI factor of value among these students.

Another study investigating the students' attitudes toward mathematics and its relationship with academic achievement in Asia has been conducted by~\cite{lim2010mathematics}. A total of 984 junior students in Singapore is measured in the study. Although they had positive attitudes toward mathematics,
it is observed that they lacked intrinsic motivation to do the subject. The author also discovered that there  was
a significant positive correlation between intrinsic motivation and academic achievement in mathematics.
It is noted that more parents and educators should focus to build intrinsic motivation among the students in order boost high self-confidence which in turn leads to greater academic achievement. Since the ATMI developed by~\cite{tapia2004instrument} is too long, a shortened version of the ATMI that only measures four subscales has been developed by~\cite{lim2013development}. The authors focused on the enjoyment of mathematics, motivation to do mathematics, self-confidence in mathematics and perceived value of mathematics. The study was also conducted in Singapore with more than 1600 participants. The result suggests that a very high correlation ($r = 0.96$) was discovered between the enjoyment and the motivation subscales. There was also a strong correlation between this shortened version of ATMI and the original version of ATMI, with a correlation mean of 0.96.

A quasi-experimental study on whether participation in single-sex classrooms versus mixed-sex classrooms makes a difference 
in attitudes toward mathematics and academic performance in PreCalculus in an American public school has been investigated by~\cite{blechle2007attitudes}.
The result indicates that there are no statistically significant differences in both the attitude toward mathematics as well as in their PreCalculus achievement.
Positive and negative attitudes toward mathematics have been redefined by~\cite{zan2008attitude}.
The authors also presented some results from the Italian project on the negative attitude toward mathematics in the context of the dichotomy between positive and negative attitudes.
A study on the possible formation of an attitude toward mathematics and the factors that might contribute to this has been discussed by~\cite{chowdhury2012identification}. Apart from gender difference and students' academic achievement, the authors also consider other factors,
including home environment, each of the parents' attitude, the language medium of instruction and the class environment.

An effort to impact student attitude toward mathematics through problem-based learning (PBL) has been attempted by~\cite{wade2013impacting}.
The author discovered that although PBL did not improve the majority of students' attitude toward mathematics,
it has a positive impact to help the student to see mathematics as something useful in the real world.
Student indicators for success in the entry-level college Calculus have been examined by~\cite{pyzdrowski2013readiness}.
From the authors' analysis, it is discovered that high school academic performance, the Calculus Readiness Assessment and the ATMI score
had positive significant correlations with course performance. 
A solid and robust instrument to measure attitude toward mathematics with contrasted evidence of validity and reliability is
presented by~\cite{palacios2013attitudes}.

A relationship between mathematics anxiety and attitude towards mathematics has been studied by~\cite{jackson2012math}.
The author focused her study among pre-service elementary teachers and investigated whether there exists a possibility to reduce 
mathematics anxiety and to improve the attitude toward mathematics.
The relationships among mathematical anxiety, attitude toward learning mathematics and separate and connected ways of knowing has been investigated 
in a doctoral dissertation by~\cite{burnes2014understanding}.
The correlation between mathematics anxiety, attitudes toward mathematics, and mathematics achievements has been investigated in a Master thesis by~\cite{chong2014adolescents}. The study found that adolescents from six secondary schools in Sarawak, East Malaysia, faced a moderate level of mathematics anxiety, have positive attitudes toward mathematics and perceived themselves to have a moderate level of achievement in mathematics.

An effort to improve students' attitude toward mathematics can be conducted in many different ways.
For instance, by incorporating a teacher's diary into the instruction portfolio~\cite{cutler1999you},
by using a computer-based instructional mathematics simulation game~\cite{van2000effect,touparova2000teaching}.
A recent proposal is regarding situated, authentic problem solving, a model that explains how digital games can promote transfer and improve attitude toward mathematics~\cite{van2015saps}. 
A comparison of mathematics achievement of 34 students in fifth grade using traditional textbook instruction and differentiated instruction 
as well as whether there exists a difference in terms of students' attitude toward mathematics after the implementation of these two different instructions 
has been studied by~\cite{gamble2011impact}. The author observed a significant difference in students' attitudes toward mathematics relative to the enjoyment component. An effort to create a feminist mathematics environment in order to engage more female students to have a more positive attitude toward mathematics has been 
attempted by~\cite{burroughs2015diversity}. 

Although the research presented in this paper does not really contribute to improving students' attitude toward mathematics,
it is still interesting to see some characteristics for a particular set of students where this type of research is rarely conducted.
The paper is organized as follows. Section~\ref{hypothesis} proposes a number of hypotheses related to our research. 
One hypothesis is regarding the correlation between previous academic performance in mathematics and the current students' attitude toward it.
The other two hypotheses investigate the students' attitude toward mathematics in terms of gender difference and the specializations that they enroll during the Foundation at NUFYP. Section~\ref{method} describes the method conducted in this research. The participants of this study are discussed as well as the instrument design and measurement. Section~\ref{result} presents the result of our analysis. The hypotheses proposed earlier are either confirmed or rejected.
Finally, Section~\ref{conclusion} provides conclusion and remark to our finding and observation.

\section{Hypotheses} \label{hypothesis}

The following hypotheses are proposed in this study.

\textbf{Hypothesis 1.} \textsl{There is a significant positive correlation between academic performance in mathematics and attitude toward it.}
Generally, a good, positive attitude toward a particular subject will influence an academic achievement in that particular subject and there is no exception with achievement in mathematics. The attitude toward mathematics and the academic achievement in it at the secondary level students have been reported amongst others by~\cite{ma1997reciprocal,ma2004determining,mubeen2013attitude}. In particular, \cite{ma2004determining} apply a statistical method known as structural equation modeling and they conclude that prior achievement of students significantly predicted later attitude across Grades 7--12.
By contrast, the prior attitude toward mathematics did not meaningfully predict later achievement.

\textbf{Hypothesis 2.} \textsl{There is a significant difference between male and female students in terms of their attitude toward mathematics.}
There exists certain preconception that mathematics and mathematics related majors are only for male while female students are more suitable to study
Arts, Languages, Humanities, Social Science or any other subjects that do not require a higher level of mathematics. A number of studies have confirmed this misconception and data on student body across many universities throughout the world support further this common opinion~\cite{nosek2002math}.
In the context of ATMI, there was an overall significant effect of gender on two of the factors of ATMI. It is observed that male students scored higher than female students on self-confidence and value~\cite{tapia2000effect,tapia2001effect}. 

A study of Ghanian secondary school students on gender difference revealed that there was a significant difference in attitudes shown toward mathematics among male and female students, where the former ones have higher ATMI scores than the latter ones~\cite{asante2012secondary}.
Another study on college male and female students on mathematics achievement and attitudes toward Remedial Mathematics in a metropolitan community college setting in the Southeastern United States has been analyzed by~\cite{hughes2016differences}.
The study found that there was a statistically significant difference between male and female students in Remedial Mathematics achievement, being male students have higher scores than female students. A similar characteristic is also observed for the attitude toward mathematics.
On the other hand, it is discovered that there is no gender difference when it comes to attitude toward mathematics among several secondary schools in Maldives~\cite{mohamed2011secondary} and in Portugal~\cite{mata2012attitudes}, although it is observed that female students exhibited a continuous decline in attitude the further they progressed in school. 

\textbf{Hypothesis 3.} \textsl{There is a significant difference among the students' specializations in terms of their attitude toward mathematics.}
It is generally perceived that the students who specialize in the Mathematical--Physics (MP) track tend to have a better attitude toward mathematics than their peers who specialize in the Biology--Chemistry (BC) track or the International Relation--Economics (IRE) track. The students enrolled in the MP track take one mathematics module that covers PreCalculus (quadratic equations, polynomials, trigonometry) and single-variable Calculus (limit, derivative, integral, introduction to differential equations). The students enrolled in the BC track do not take any mathematics module at all, while the students enrolled in the IRE track take one mathematics module called `Mathematics for Economics'. The content is similar to the module for the MP track, but this module excludes an introduction to differential equations. Rather, it includes more in-depth coverage of Statistics, which covers descriptive statistics, basic probability, permutations and combinations, binomial and Poisson distributions.

To the best of our knowledge, there is no literature discusses these particular three specializations since the program is unique only to a setup university in Kazakhstan, for which the Foundation program is supported by the UCL and thus adopting a British style curriculum. There exists, however, a study on how family socioeconomic status and parental involvement affect college major choices among high school students in the United States~\cite{ma2009family}. The students' attitude toward different subjects is measured in the tenth grade while the students' achievement in Mathematics and English are measured in the twelfth grade. The author discovered that parental involvement in children's domain-specific education exerts significant effects on children's college major choice.

\section{Method} \label{method}

\subsection{Participant}

The participants of this study are the students enrolled in the NUFYP in Astana, Kazakhstan, during the academic year of 2012/2013. 
This means that these students are the cohort of the graduating Class 2017 from the undergraduate program at the same university. 
Depending on the age they start elementary school, which typically around age six or seven, they completed the secondary school of Grade~11 at the age of 17 or 18. Hence, at the end of their Foundation period, the age of participants is between 18 and 19 years old.
An estimated number of targets are around 500 students. A sample of 108 valid questionnaires was received.
According to gender groups, there are 55 males and 53 females. According to their specialization, 13 students were enrolled in IRE track, 22 students were enrolled in the BC track and 73 students were enrolled in the MP track. Out of 108 students, only seven students received Grade~4 (the second highest) for their Mathematics score in high school, while the rest (i.e. 101 students), received Grade~5 (the highest). The composition of polled students according to gender and specializations is shown in Table~\ref{polled}.
\begin{table}[h]
\begin{center}
\begin{tabularx}{0.6\textwidth}{@{}X*{4}{R}@{}}
\toprule
       & IRE     & BC       & MP        & Total \\ \hline
Male   &  2      &  6       & 47        & 55    \\
Female & 11      & 16       & 26        & 53    \\ 
Total  & 13      & 22       & 73        & 108   \\ \bottomrule
\end{tabularx}
\end{center} 
\caption{The composition of polled students according to gender and specializations.} \label{polled}
\end{table}

\subsection{Instrument design and measurement}

An online questionnaire designed using a Google document was distributed to a common university, 
electronic mail address to the students enrolled in NUFYP on 17 June 2013.
The questionnaire remained open for a couple of days, but the final response was obtained on 18 June 2013.
The preparation of this online questionnaire was conducted by research students Damira Alshimbayeva and Korkem Ybrakhym during their summer internship at our school SST. Both are Mathematics major students from the graduating Class of 2016.
Since usually the intake in the NUFYP is around 500, it is estimated that the total target of students is around 500.
The questionnaire collects basic information about gender, study program specialization (track) and secondary school mathematics score.
The later one is ordered from Grade~1 (the lowest score) to Grade~5 (the highest score).
The rest of the questionnaire consists of 40 ATMI items with a five-item Likert scale.
The options~1~until~5 indicate `strongly disagree' to `strongly agree', respectively. Option~3 means `neutral'.
The ATMI utilized in the data collection is presented in the Appendix section.

\subsection{Statistical analysis}
A number of simple statistical analyses were conducted in this study.
This includes calculating mean, standard deviation, descriptive Pearson correlation, $t$-test and analysis of variance (ANOVA).
Data analysis was performed using \textsl{LibreOffice Calc}, \textsl{Geogebra} and a statistical computing software \textsl{R}.
The histograms are produced using \textsl{R} too.

\section{Result} \label{result}

\subsection{ATMI score}

\begin{figure}[h]
\begin{center}
\includegraphics[width = 0.85\textwidth]{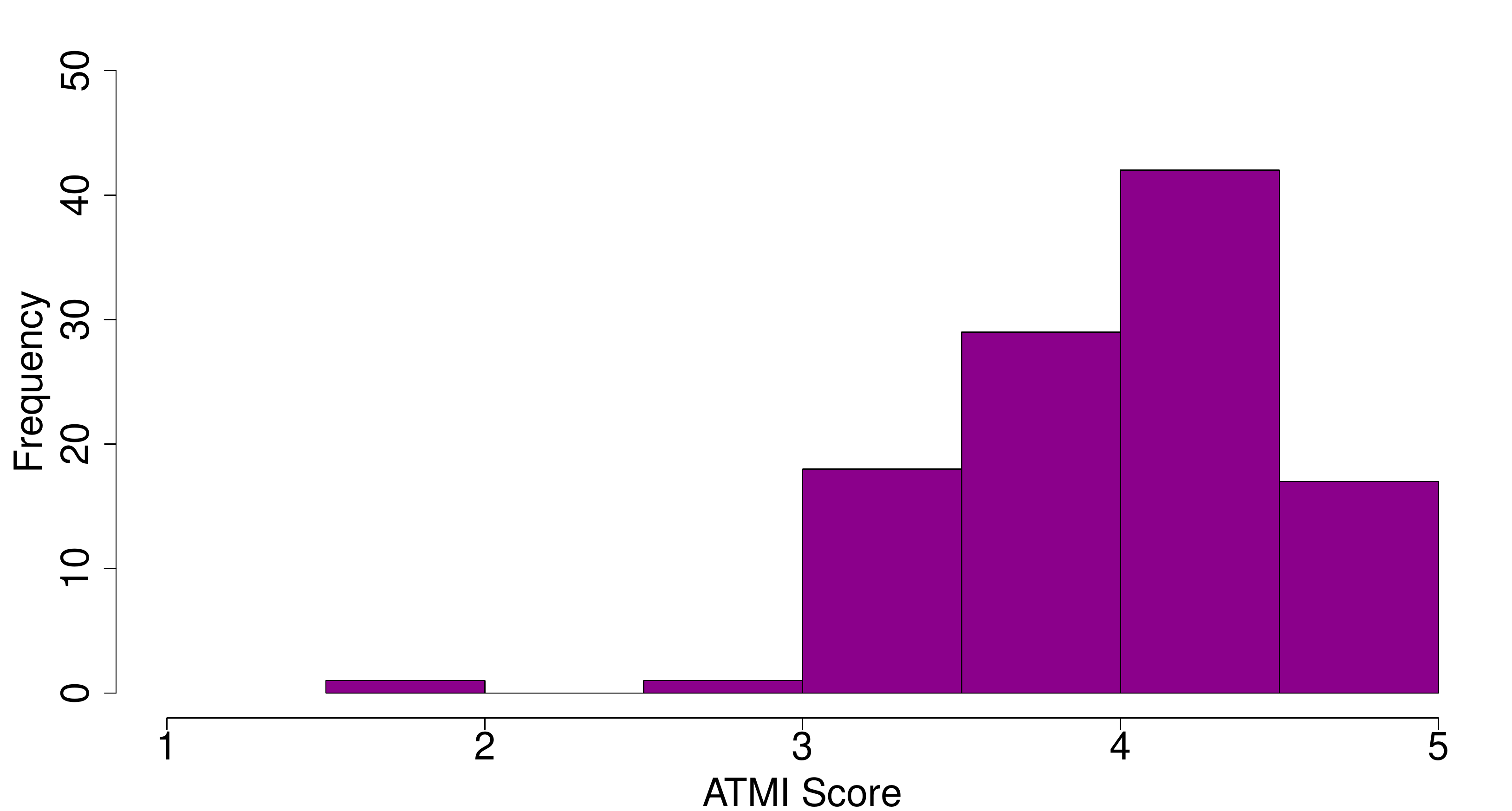}
\end{center}
\caption{Histogram of ATMI scores for all students.} \label{atmiall}
\end{figure}

The ATMI score indicates the students' attitude toward mathematics.
The minimum and the maximum ATMI scores are one and five, respectively.
Thus, unless indicated specifically, in the context of this article, a higher or a lower ATMI score indicates a more positive or more negative attitude toward mathematics and vice versa. 
From the obtained polled students, the mean of the ATMI score is 3.999/5.00, with its standard deviation of 0.531.
Hence, the students enrolled at NUFYP generally have a positive attitude toward mathematics, as shown by relatively high scores of ATMI.
The distribution of ATMI scores for all students, male, female, the students who specialized in the MP, BC and IRE tracks are presented in Table~\ref{alldata}.
More than 98\% of the students show ATMI scores of three and above and among these, more than 50\% display the ATMI scores of between four and five.
The histogram of the ATMI scores for all polled students is presented in Figure~\ref{atmiall}.
\begin{table}[h]
\begin{center}
\begin{tabularx}{0.8\textwidth}{@{}c*{5}{R}c@{}}
\toprule
ATMI score & All & Male & Female & MP & BC & IRE \\ \hline
$(1,1.5]$  & 0   & 0    & 0      & 0  & 0  & 0   \\
$(1.5,2]$  & 1   & 0    & 1      & 0  & 0  & 1   \\
$(2,2.5]$  & 0   & 0    & 0      & 0  & 0  & 0   \\
$(2.5,3]$  & 1   & 1    & 0      & 1  & 0  & 0   \\
$(3,3.5]$  & 18  & 8    & 10     & 10 &  4 & 4   \\
$(3.5,4]$  & 29  & 15   & 14     & 15 & 11 & 3   \\
$(4,4.5]$  & 42  & 21   & 21     & 32 &  7 & 3   \\
$(4.5,5]$  & 17  & 10   & 7      & 15 &  0 & 2   \\
\bottomrule
\end{tabularx}
\end{center}
\caption{The ATMI score distribution for all students, male, female, the students who specialized in MP, BC, and IRE tracks.} \label{alldata}
\end{table}

\subsection{Academic performance}

In this subsection, Hypothesis~1 (there is a significant positive correlation between academic performance in mathematics and attitude toward it) is examined.
What it means by the term `academic performance' is the previous academic performance of the students during their high school year.
The data only collected their final grade of mathematics at the secondary school level. 
It turns out that these data are not particularly useful in terms of score diversity since the majority of the students enrolled in NUFYP are among the best in the country and thus it generally translates to excellent score in mathematics during high school.

Out of 108 valid returned questionnaires, 101 of the students obtained the highest score (Grade~5) and only seven students have the second highest score (Grade~4).
The average of high school mathematics grade is 4.935/5.00.
After calculating the ATMI score of these students, in general, they have a positive attitude toward mathematics.
The mean of the ATMI score is 3.999 on the scale from one to five, with its standard deviation of 0.531.
Further analysis confirms Hypothesis~1. There is indeed a positive correlation between previous high achievement in mathematics and a favorable attitude toward mathematics. The Pearson correlation value is $r = 0.300$ and although its effect size is medium, it is extremely significant ($p$-value $= 0.0008 < 0.001$).
Thus, since the result is extremely significant, Hypothesis~1 is accepted.

\subsection{Gender difference}

\begin{figure}[h]
\begin{center}
\includegraphics[width = 0.75\textwidth]{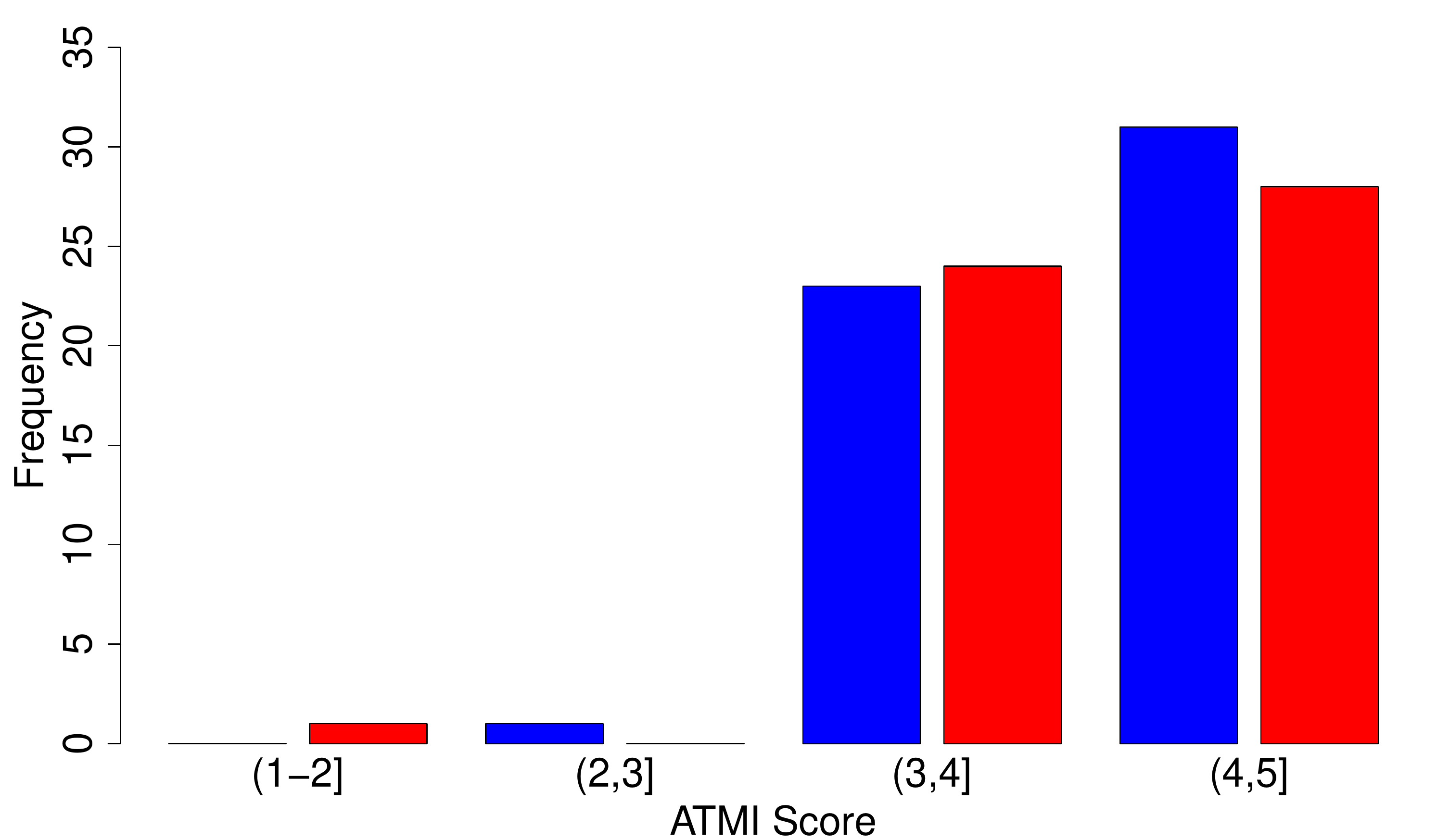}
\end{center}
\caption{(Color online) Histogram of ATMI scores for male (blue, left-hand side bars) and female (red, right-hand side bars) students.} \label{gender01}
\end{figure}

In this subsection, Hypothesis~2 (there is a significant difference between male and female students in terms of their attitude toward mathematics) is examined.
Out of 55 male students and 53 female students, the ATMI scores are 4.0332 and 3.9637 out of five, respectively.
The histogram of ATMI scores for male and female students is displayed in Figure~\ref{gender01}.
From this result, it seems that male students have a better attitude toward mathematics than their female counterparts.
However, the unpooled $t$-test gives the standard error $SE = 0.1023$, degrees of freedom $= 104.6697$, $t = 0.6795$ and $p$-value $= 0.4983 > 0.05$ (not significant).
Since the result is not significant, there is no sufficient evidence to accept the hypothesis.
Thus, when it comes to gender difference, there is no significant difference between male and female students when it comes to their attitude toward mathematics.
The summary of the unpooled $t$-test is shown in Table~\ref{ttest}. Since the result is not significant, Hypothesis~2 is rejected.
\begin{table}[h]
\begin{center}
\begin{tabularx}{0.8\textwidth}{@{}l*{2}{c}R@{}}
\toprule
Statistic    	     		& Male 		& Female & \\ \hline
Sample mean 				& 4.0332 	& 3.9637 & \\ 
Sample standard deviation \qquad \qquad & 0.5106    & 0.5507 & \\
Number of data				& 55        & 53     & \\ 
Standard error (SE)         & 		    &        & 0.1023 \\
Degrees of freedom ($df$)   & 		    &        & 104.6697 \\
$t$-Value   				& 		    &	     & 0.6795 \\
$p$-Value					& 		    &        & 0.4983 \\ 
\bottomrule
\end{tabularx}
\end{center}
\caption{A summary of $t$-test by comparing the mean difference between male and female students.
It gives the values of $t = 0.6795$ and $p = 0.4983 > 0.05$ (not significant).} \label{ttest}
\end{table}

\subsection{Specialization}

\begin{figure}[h]
\begin{center}
\includegraphics[width = 0.75\textwidth]{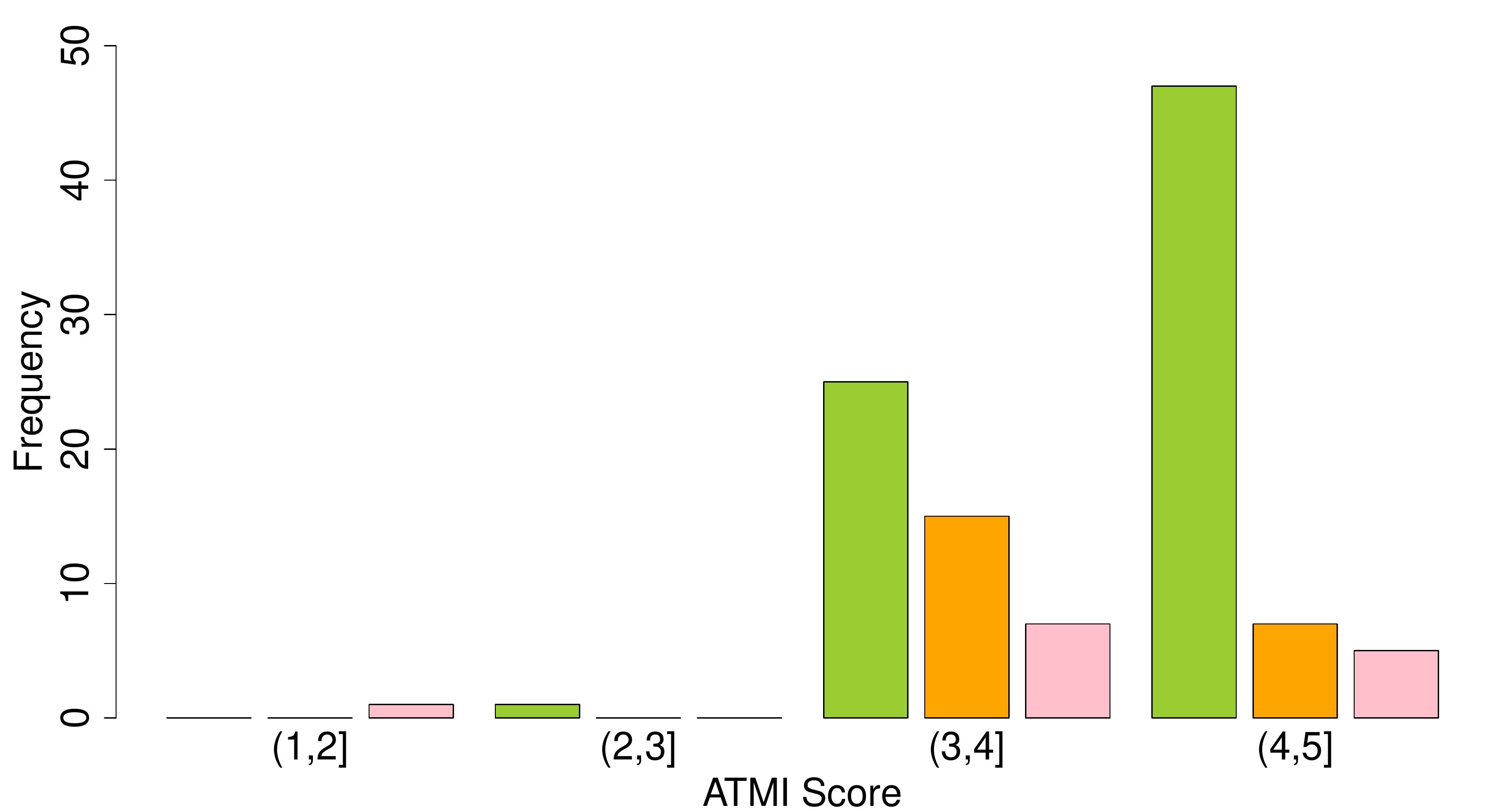}
\end{center}
\caption{(Color online) Histogram of ATMI scores for students specialized in the MP track (green, left-hand side bars), the BC track (orange, middle bars) and the IRE track (pink, right-hand side bars).} \label{track1}
\end{figure}

In this subsection, Hypothesis~3 (there is a significant difference among the students' specializations in terms of their attitude toward mathematics) is examined.
As mentioned previously, there are three different study programs, also known as specializations or tracks, that the students in NUFYP enrolled during their Foundation study, i.e. IRE track, BC track and MP track.
From the collected sample data, there are 13 students specialized in the IRE track, 22 students specialized in the BC track and 73 students specialized in the MP track. The average ATMI scores for each specialization is 3.6808, 3.8147 and 4.1113 on the scale from one to five, respectively.
It seems that the students in the MP track have the best attitude toward mathematics among the three groups,
while the students in the IRE track show the least favorable attitude toward mathematics.
The students in the BC track fall somewhere in between.
The histogram of ATMI scores for students specialized in the MP, BC and IRE tracks is displayed in Figure~\ref{track1}.

A one-way ANOVA for independent samples was conducted to analyze these data and the standard weighted mean analysis is adopted.
The result shows that there is no significant difference in ATMI scores between the students in the IRE track and the ones in the BC track,
as well as between the students in the BC track and the ones in the MP track.
However, there is a significant difference in ATMI scores between the students in the IRE track and the ones in the MP track.
A further analysis is conducted using the Tukey's honestly significant difference (HSD) test.
The result $F(2,105) = 5.6848$ with $p$-value $= 0.0045 < 0.01$ are obtained.
The summary of the ANOVA and post-hoc Tukey's HSD test is shown in Table~\ref{anova}.
Since there exists a very significant result, Hypothesis~3 is accepted.
\begin{table}[h]
\begin{center}
\begin{tabularx}{0.8\textwidth}{@{}l*{2}{r}*{3}c@{}}
\toprule
Source     & $df${\:} & SS {\;\;}& MS     & $F$    & $p$-Value \\ \hline
Treatments & 2    & 2.9838   & 1.4919 & 5.6848 & 0.0045 \\
Error      & 105  & 27.5561  & 0.2624 &        &        \\
Total      & 107  & 30.5399  &        &        &        \\
\bottomrule
\end{tabularx}
\end{center}
\caption{ANOVA summary with three independent samples (IRE, BC and MP tracks).
The post-hoc Tukey's HSD test gives $F(2,105) = 5.6848$ and $p$-value $= 0.0045 < 0.01$ (very significant).} \label{anova}
\end{table}

This result is somewhat surprising and rather unexpected. Since only two tracks (MP and IRE) offer modules in mathematics, 
we would expect that the students in the BC track have the least favorable attitude toward mathematics.
We also would expect that there is a significant difference when it comes to the ATMI score between the students enrolled in the MP and the BC tracks as well as 
the students enrolled in the BC and the IRE tracks.
But it turns out that there are no significant differences between these two respective groups in terms of their attitude toward mathematics. 
A possible explanation for this finding is that the number of students enrolled in the MP track 
who participated in this study is more than three times than the number of students enrolled the BC track (73\,:\,22) and 
more than five times than the number of students enrolled in the IRE track (73\,:\,13).
Therefore, any exhibit of a more favorable attitude toward mathematics among the students in the IRE track is shielded by the students in the MP track.

\section{Conclusion} \label{conclusion}

The attitude toward mathematics among the students enrolled at NUFYP has been investigated in this paper.
In the following, we will discuss some weaknesses in our research, draw a conclusion, provide remarks and possible future direction of the research.
One weakness of this research is the relatively low response rate to the online questionnaire among the students enrolled at NUFYP.
From around 500 students targeted in this research, only 108 valid responses are returned. Thus, the response rate is only 21.6\%.
One possible explanation is the timing of the data collection. The month of June is near the end of the study period for the NUFYP students.
Therefore, the majority of the students, more than three-quarter of them, prefer to ignore the questionnaire and focus on their study instead.

The second weakness of this research is the inquiry of the previous academic achievement in mathematics, 
namely the collection of the final mathematics score during the secondary school period.
It is discovered that there is a very little variation in the score distribution among the polled students.
The majority of the students obtained the highest grade (Grade 5) and a small number obtained the second highest grade (Grade 4).
The correlation between the academic achievement in mathematics and the attitude toward it would be having a better visibility if 
the mathematics achievement of the students was collected after they enrolled at NUFYP.
Collecting a diverse mathematics score, let say in the range of zero to 100, will provide more variation rather than the high school score in the range from one to five with only integer scores.

Another weakness lies in the distribution of students in terms of their specializations. 
The collected data of students' specializations capture more than two third (67.6\% or 73 students) of the students specialized in the MP track,
while only 20\% (22 students) and 12\% (13 students) of the students specialized in the BC and IRE tracks, respectively.
There is a chance that our statistical result provides a different outcome if there are more or less equal distribution of specializations among the polled students.
Apart from academic achievement, gender difference, and study program specialization, it would be interesting to see other backgrounds too, including the hometown or the place of the origin (urban or rural) and socioeconomic status of the students. Compare with~\cite{ma2009family,chowdhury2012identification}.

From the statistical analysis that we have conducted, we can draw the following conclusion.
The result shows that there is a positive correlation with medium effect size between the high school mathematics score and 
the attitude toward mathematics as indicated by the ATMI score.
A statistical analysis using unpooled $t$-test suggests that there is no significant difference between gender when it comes to the attitude toward mathematics
($t = 0.6795$, $p = 0.4983 > 0.05$), even though it appears that the male students have a slightly higher ATMI score than their female counterparts.
Regarding the specialization, the students in the MP track have the highest ATMI score and the students in the IRE track has the lowest ATMI score. The students in the BC track has ATMI score between the other two groups, even though they do not take any mathematics module during the entire program.
Further analysis using one-way ANOVA for independent samples and post-hoc Tukey's HSD discovered that there is a very significant difference between the students specialized in the IRE track and the ones in the MP track when it comes to the attitude toward mathematics ($F(2,105) = 5.6848$, $p = 0.0045 < 0.01$). 
There is no significant difference between the students in the IRE and the BC tracks as well as between the students in the BC and the MP tracks.

A possible future direction of this research is to conduct a similar test to the same group of students
and investigate whether there exists a significant difference in the attitude toward mathematics between the time when they were enrolled at NUFYP and the time of nearly the end of their undergraduate study period.
Currently, these students are the graduating Class of 2017 distributed at various undergraduate programs at NU.
The previously NUFYP students will be affiliated with one of the schools at NU: SST, SEng, or SHSS.
Furthermore, they have declared their major too. Except for International Relations and Political Science as well as World Language and Literature majors,
all other majors require the students to take several mathematics courses, including but not limited to, PreCalculus, Calculus 1 and 2 and Linear Algebra with Applications.
Thus, it would be interesting to investigate whether there is a significant difference in the attitude toward mathematics 
in terms of different schools, various majors and the number of credit hours in mathematics that the students have taken.
Our conjecture is that the students enrolled in schools and majors that use a lot of mathematics will have a better attitude toward mathematics than the one with less exposure to mathematics. A similar conjecture can be proposed to the students who take more courses in mathematics will generally have a better attitude toward mathematics.

\section*{\large Acknowledgements}
{\small The author would like to thank Damira Alshimbayeva and Korkem Ybrakhym (Mathematics major students from the graduating Class 2016) for their assistance in  collecting data during the Summer internship period in 2013, as well as to Gulzhan Tumenbayeva (Mathematics major graduating Class 20017), Binur Yermukanova and Aigerim Bulambayeva (Economics major graduating Class 2017) for information about NUFYP. The students enrolled in NUFYP during the academic year 2012/2013, who are also the Nazarbayev University students from the graduating Class 2017 who have been participating in filling in the questionnaire are gratefully acknowledged.
The author would like to acknowledge Su Ting Yong (The University of Nottingham Malaysia Campus) for an inspiration to conduct research on attitude toward mathematics and Martha Tapia (Berry College, Georgia, USA) for the permission to use the ATMI in this research.
\par}

\section*{\large Disclosure statement}
{\small No potential conflict of interest was reported by the author.}

\section*{\large ORCID}
{\small {\sl N. Karjanto} \; \url{https://orcid.org/0000-0002-6859-447X}}

{\footnotesize
\bibliography{arXivAtmi} 
\bibliographystyle{apalike}
}

\vspace{1cm}
\newpage
\begin{center}
\textbf{Appendix. {\slshape Attitudes Toward Mathematics Inventory}}
\end{center}
{\small 
\parindent=0pt
\textbf{Instruction:} This inventory consists of statements about your attitude toward mathematics. There
are no correct or incorrect responses. Read each item carefully. Please think about how you feel
about each item. Enter the letter that most closely corresponds to how each statement best
describes your feelings. Please answer every question.
Please use these response codes: 1 -- Strongly Disagree, 2 -- Disagree, 3 -- Neutral, 4 -- Agree, 5 -- Strongly Agree.

\renewcommand{\tabcolsep}{6pt}
\renewcommand{\arraystretch}{1.05}
\begin{center}
\begin{tabularx}{\textwidth}{@{}cp{0.72\textwidth}*{5}{c}@{}}
\toprule
No & Statement                                               & 1 & 2 & 3 & 4 & 5 \\ \hline 
1  & Mathematics is a very worthwhile and necessary subject. & $\square$ & $\square$ & $\square$ & $\square$ & $\square$ \\ 
2  & I want to develop my mathematical skills.  & $\square$ & $\square$ & $\square$ & $\square$ & $\square$ \\ 
3  & I get a great deal of satisfaction out of solving a mathematics problem. & $\square$ & $\square$ & $\square$ & $\square$ & $\square$ \\ 
4  & Mathematics helps develop the mind and teaches a person to think. & $\square$ & $\square$ & $\square$ & $\square$ & $\square$ \\ 
5  & Mathematics is important in everyday life. & $\square$ & $\square$ & $\square$ & $\square$ & $\square$ \\ 
6  & Mathematics is one of the most important subjects for people to study. & $\square$ & $\square$ & $\square$ & $\square$ & $\square$ \\ 
7  & High school math courses would be very helpful no matter what \break I decide to study. & $\square$ & $\square$ & $\square$ & $\square$ & $\square$ \\ 
8  & I can think of many ways that I use math outside of school. & $\square$ & $\square$ & $\square$ & $\square$ & $\square$ \\ 
9  & Mathematics is one of my most dreaded subjects. & $\square$ & $\square$ & $\square$ & $\square$ & $\square$ \\ 
10 & My mind goes blank and I am unable to think clearly when working with \break mathematics. & $\square$ & $\square$ & $\square$ & $\square$ & $\square$ \\ 
11 & Studying mathematics makes me feel nervous. & $\square$ & $\square$ & $\square$ & $\square$ & $\square$ \\ 
12 & Mathematics makes me feel uncomfortable. & $\square$ & $\square$ & $\square$ & $\square$ & $\square$ \\ 
13 & I am always under a terrible strain in a math class. & $\square$ & $\square$ & $\square$ & $\square$ & $\square$ \\ 
14 & When I hear the word mathematics, I have a feeling of dislike. & $\square$ & $\square$ & $\square$ & $\square$ & $\square$ \\ 
15 & It makes me nervous to even think about having to do a mathematics \break problem. & $\square$ & $\square$ & $\square$ & $\square$ & $\square$ \\ 
16 & Mathematics does not scare me at all. & $\square$ & $\square$ & $\square$ & $\square$ & $\square$ \\ 
17 & I have a lot of self-confidence when it comes to mathematics. & $\square$ & $\square$ & $\square$ & $\square$ & $\square$ \\ 
18 & I am able to solve mathematics problems without too much difficulty. & $\square$ & $\square$ & $\square$ & $\square$ & $\square$ \\ 
19 & I expect to do fairly well in any math class I take. & $\square$ & $\square$ & $\square$ & $\square$ & $\square$ \\ 
20 & I am always confused in my mathematics class. & $\square$ & $\square$ & $\square$ & $\square$ & $\square$ \\ 
21 & I feel a sense of insecurity when attempting mathematics. & $\square$ & $\square$ & $\square$ & $\square$ & $\square$ \\ 
22 & I learn mathematics easily. & $\square$ & $\square$ & $\square$ & $\square$ & $\square$ \\ 
23 & I am confident that I could learn advanced mathematics. & $\square$ & $\square$ & $\square$ & $\square$ & $\square$ \\ 
24 & I have usually enjoyed studying mathematics in school. & $\square$ & $\square$ & $\square$ & $\square$ & $\square$ \\ 
25 & Mathematics is dull and boring. & $\square$ & $\square$ & $\square$ & $\square$ & $\square$ \\ 
26 & I like to solve new problems in mathematics. & $\square$ & $\square$ & $\square$ & $\square$ & $\square$ \\ 
27 & I would prefer to do an assignment in math than to write an essay. & $\square$ & $\square$ & $\square$ & $\square$ & $\square$ \\ 
28 & I would like to avoid using mathematics in college. & $\square$ & $\square$ & $\square$ & $\square$ & $\square$ \\
29 & I really like mathematics. & $\square$ & $\square$ & $\square$ & $\square$ & $\square$ \\ 
30 & I am happier in a math class than in any other class. & $\square$ & $\square$ & $\square$ & $\square$ & $\square$ \\ 
31 & Mathematics is a very interesting subject. & $\square$ & $\square$ & $\square$ & $\square$ & $\square$ \\ 
32 & I am willing to take more than the required amount of mathematics. & $\square$ & $\square$ & $\square$ & $\square$ & $\square$ \\ 
33 & I plan to take as much mathematics as I can during my education. & $\square$ & $\square$ & $\square$ & $\square$ & $\square$ \\ 
34 & The challenge of math appeals to me. & $\square$ & $\square$ & $\square$ & $\square$ & $\square$ \\ 
35 & I think studying advanced mathematics is useful. & $\square$ & $\square$ & $\square$ & $\square$ & $\square$ \\ 
36 & I believe studying math helps me with problem solving in other areas. & $\square$ & $\square$ & $\square$ & $\square$ & $\square$ \\ 
37 & I am comfortable expressing my own ideas on how to look for solutions \break to a difficult problem in math. & $\square$ & $\square$ & $\square$ & $\square$ & $\square$ \\ 
38 & I am comfortable answering questions in math class. & $\square$ & $\square$ & $\square$ & $\square$ & $\square$ \\ 
39 & A strong math background could help me in my professional life. & $\square$ & $\square$ & $\square$ & $\square$ & $\square$ \\ 
40 & I believe I am good at solving math problems. & $\square$ & $\square$ & $\square$ & $\square$ & $\square$ \\ 
\bottomrule
\end{tabularx}
\end{center}
}
\end{document}